%% file: main.tex
\documentclass[12pt,a4paper,oneside]{article}

\usepackage{styling}

\begin{document}

\title{Unipotent elements forcing irreducibility in linear algebraic groups}
\author{Mikko Korhonen\thanks{Section de mathématiques, École Polytechnique Fédérale de Lausanne, CH-1015 Lausanne, Switzerland} \thanks{Email address: mikko.korhonen@epfl.ch}}
\date{\today}

\maketitle

\begin{abstract}Let $G$ be a simple algebraic group over an algebraically closed field $K$ of characteristic $p > 0$. We consider connected reductive subgroups $X$ of $G$ that contain a given distinguished unipotent element $u$ of $G$. A result of Testerman and Zalesski (Proc. Amer. Math. Soc., 2013) shows that if $u$ is a regular unipotent element, then $X$ cannot be contained in a proper parabolic subgroup of $G$. We generalize their result and show that if $u$ has order $p$, then except for two known examples which occur in the case $(G, p) = (C_2, 2)$, the subgroup $X$ cannot be contained in a proper parabolic subgroup of $G$. In the case where $u$ has order $> p$, we also present further examples arising from indecomposable tilting modules with quasi-minuscule highest weight.\end{abstract}

\section{Introduction}
\input{PROJECT_intro}

\section{Notation}\label{section:notation}
\input{PROJECT_preliminaries}

\section{Tilting modules with quasi-minuscule highest weight}\label{section:quasiminuscule}
\input{PROJECT_quasiminuscule}

\section{Representations of $\SL_2(K)$}\label{section:SL2section}
\input{PROJECT_SL2}

\section{Distinguished unipotent elements of order $p$ in non-$G$-cr subgroups}\label{section:nonGcr_orderp}
\input{PROJECT_nonGcr}

\section{Further examples}\label{section:examples_section}
\input{PROJECT_examples_v2}

\section*{Acknowledgements}
\input{PROJECT_acknowledgements}

\bibliographystyle{alpha}
\bibliography{bibliography} 

\end{document}

%% file: PROJECT_intro.tex
%

Let $G$ be a simple linear algebraic group over an algebraically closed field $K$ of characteristic $p > 0$. A unipotent element $u \in G$ is said to be \emph{distinguished} if its centralizer does not contain a nontrivial torus. We say that $u$ is \emph{regular}, if the dimension of its centralizer is equal to the rank of $G$. It is well known that regular unipotent elements are distinguished \cite[1.14(a)]{SpringerSteinberg}.

One approach towards understanding the subgroup structure and the properties of the unipotent elements of $G$ is to study the subgroups which contain a fixed unipotent element $u \in G$. There are many results in this direction in the literature, which have found application in linear algebraic groups and elsewhere. To give one example, Saxl and Seitz \cite{SaxlSeitz} classified the positive-dimensional maximal closed overgroups of regular unipotent elements; this classification and further study of overgroups of regular unipotent elements was applied to the inverse Galois problem by Guralnick and Malle in \cite{GuralnickMalle}.

This paper is concerned with the overgroups of distinguished unipotent elements and is part of ongoing work by the present author, which aims to classify all connected reductive overgroups of distinguished unipotent elements of $G$. In this paper, we consider generalizations of the following result of Testerman and Zalesski on overgroups of regular unipotent elements to other classes of distinguished unipotent elements. Below a subgroup of $G$ is said to be \emph{$G$-irreducible} ($G$-ir), if it is not contained in any proper parabolic subgroup of $G$.

\begin{lause}[{\cite{TestermanZalesski}}]\label{thm:testerman_zalesski}
Let $X$ be a connected reductive subgroup of the simple algebraic group $G$. If $X$ contains a regular unipotent element of $G$, then $X$ is $G$-irreducible.
\end{lause}

As an application Theorem \ref{thm:testerman_zalesski} and the Saxl-Seitz classification of maximal closed reductive subgroups of $G$ that contain a regular unipotent element \cite[Theorem A, Theorem B]{SaxlSeitz}, one can classify all connected reductive subgroups of $G$ that contain a regular unipotent element \cite[Theorem 1.4]{TestermanZalesski}. Thus with the classification of connected reductive overgroups of distinguished unipotent elements in mind, one might hope that Theorem \ref{thm:testerman_zalesski} would generalize to all distinguished unipotent elements, but this turns out not to be the case. The smallest examples are given by the following result.

%
%

\begin{restatable}[]{prop}{propSP}
\label{prop:resultSP4_intro}
Assume that $p = 2$. Let $G = \Sp(V)$, where $\dim V = 4$, so $G$ is simple of type $C_2$ and has a unique conjugacy class of distinguished unipotent elements of order $p$. Let $u \in G$ be a distinguished unipotent element of order $p$. If $X < G$ is connected reductive and $u \in X$, then $X$ is $G$-irreducible unless $X$ is conjugate to $X'$, where $X'$ is described by one of the following:

\begin{enumerate}[\normalfont (i)]
\item $X'$ is simple of type $A_1$ with natural module $E$, embedded into $G$ via $E \perp E$ (orthogonal direct sum).
\item $X'$ is simple of type $A_1$ with natural module $E$, embedded into $G$ via $E \otimes E$.
\end{enumerate}

Furthermore, subgroups $X'$ in (i) and (ii) exist, contain a conjugate of $u$, are contained in a proper parabolic subgroup, and the conditions (i) and (ii) determine $X'$ up to conjugacy in $G$. 
\end{restatable}

Our main result is the following, which shows that Theorem \ref{thm:testerman_zalesski} holds for distinguished unipotent elements of order $p$, except for the two examples given by Proposition \ref{prop:resultSP4_intro}.


\begin{restatable}[]{lause}{mainresult}
\label{theorem:main_result_of_paper}
Let $u \in G$ be a distinguished unipotent element of order $p$. Let $X$ be a connected reductive subgroup of $G$ containing $u$. Then $X$ is $G$-irreducible, unless $p = 2$, the group $G$ is simple of type $C_2$, and $X$ is simple of type $A_1$ as in Proposition \ref{prop:resultSP4_intro} (i) or (ii).
\end{restatable}

In the (submitted) PhD thesis of the present author, a complete list of maximal closed connected overgroups of all distinguished unipotent elements (including those of order $> p$) is given. Combining this with Theorem \ref{theorem:main_result_of_paper}, one can give a description of all connected reductive subgroups of $G$ that contain a distinguished unipotent element of order $p$. For distinguished unipotent elements of order $> p$, we will establish in the final section of this paper some further examples, which show that the statement of Theorem \ref{thm:testerman_zalesski} does not generalize to distinguished unipotent elements. However, we believe that such examples are quite rare, and a complete classification of them should eventually be possible.

\begin{remark}In characteristic zero, the problem of classifying all connected reductive overgroups of distinguished unipotent elements is solved in the following sense. The maximal \emph{connected reductive} overgroups can be found using \cite{LiebeckSeitzTesterman} and \cite{LawtherFusion}. Furthermore, it follows from a theorem of Mostow \cite{Mostow} that any connected reductive subgroup of $G$ is contained in a \emph{reductive} maximal closed connected subgroup of $G$. Therefore in characteristic zero, we have a recursive description of the connected reductive overgroups of distinguished unipotent elements.\end{remark}

The notion of irreducibility in algebraic groups given above is due to Serre. We will also need his definition of complete reducibility in algebraic groups. For an overview of these concepts and known results, see \cite{Serre1998}, \cite{Serre}, and \cite{BMSGeometric}.

\begin{maar}[{Serre, \cite{Serre1998}}]\label{def:nonGcr}
Let $H$ be a closed subgroup of $G$. We say that $H$ is \emph{$G$-completely reducible} ($G$-cr), if whenever $H$ is contained in a parabolic subgroup $P$ of $G$, it is contained in a Levi factor of $P$. Otherwise we say that $H$ is \emph{non-$G$-completely reducible} (non-$G$-cr). 
\end{maar}

We now give the outline for our proof of Theorem \ref{theorem:main_result_of_paper}. Let $X$ be a connected reductive subgroup of $G$ containing a distinguished unipotent element $u \in G$ of order $p$. 

If $p \neq 3$ or if $X$ does not have a simple factor of type $G_2$, it follows from a result of Testerman \cite{TestermanA1} and Proud-Saxl-Testerman \cite[Theorem 5.1]{ProudSaxlTesterman} that $u$ is contained in a simple subgroup $X' < X$ of type $A_1$. Once we have established Theorem \ref{theorem:main_result_of_paper} for subgroups of type $A_1$, we know that either 

\begin{itemize}
\item $X'$ is $G$-irreducible, and thus so is $X$; or
\item $X'$ is as in Proposition \ref{prop:resultSP4_intro} (i) or (ii), and in this case it is straightforward to see that either $X = X'$ or $X$ is $G$-irreducible.
\end{itemize}

The case where $p = 3$ and $X$ has a simple factor of type $G_2$ is easy to deal with, so the question is reduced to the case where $X$ is simple of type $A_1$.

When $X$ is simple of type $A_1$ and $G$ is simple of classical type, Theorem \ref{theorem:main_result_of_paper} can be reformulated as a result in the representation theory of $\SL_2(K)$. In Section \ref{section:SL2section}, Proposition \ref{fewblocksofsizep}, we classify all $\SL_2(K)$-modules on which a non-identity unipotent element $u$ of $\SL_2(K)$ acts with at most one Jordan block of size $p$. As a consequence, we find that a self-dual $\SL_2(K)$-module is semisimple if $u$ acts on it with at most one Jordan block of size $p$, allowing us to establish Theorem \ref{theorem:main_result_of_paper} when $p \neq 2$ and $G$ is simple of classical type. When $p = 2$ and $G$ is simple of classical type, it is easily seen that distinguished unipotent elements of order $p$ exist only in the case where $G$ is simple of type $C_2$. This case is treated in the proof of Proposition \ref{prop:resultSP4_intro}.

For $G$ simple of exceptional type in good characteristic, Litterick and Thomas \cite{LitterickThomasARXIV} have classified all connected reductive non-$G$-cr subgroups, in particular all non-$G$-cr subgroups $X$ of type $A_1$. They have also described the action of such $X$ on the minimal-dimensional and adjoint modules of $G$, which together with a result of Lawther \cite{Lawther} is enough information to determine the conjugacy class of a unipotent element $u \in X$ in $G$. One finds that none of the non-$G$-cr subgroups of type $A_1$ contain distinguished unipotent elements, and it is not too difficult to see that the same is true in bad characteristic as well. In other words, we find that every subgroup $X$ of type $A_1$ containing a distinguished unipotent element of order $p$ is $G$-cr. Since $u$ is distinguished, it is easy to deduce that any such subgroup $X$ is in fact $G$-irreducible (Lemma \ref{lemma:basicgcrdist}), so Theorem \ref{theorem:main_result_of_paper} follows.


%% file: PROJECT_preliminaries.tex
%

We fix the following notation and terminology. Throughout the text, let $K$ be an algebraically closed field of characteristic $p > 0$. All the groups that we consider are linear algebraic groups over $K$, and $G$ will always denote a simple linear algebraic group over $K$. By a subgroup we will always mean a closed subgroup, and by a $G$-module $V$ we will always mean a finite-dimensional rational $KG$-module. We say that $p$ is a \emph{good prime for $G$}, if $G$ is simple of type $B_l$, $C_l$, or $D_l$, and $p > 2$; if $G$ is simple of type $G_2$, $F_4$, $E_6$, or $E_7$, and $p > 3$; or if $G$ is simple of type $E_8$ and $p > 5$. Otherwise we say that $p$ is a \emph{bad prime for $G$}.

We fix a maximal torus $T$ of $G$ with character group $X(T)$. Fix a base $\Delta = \{ \alpha_1, \ldots, \alpha_l \}$ for the roots of $G$, where $l$ is the rank of $G$. Here we use the standard Bourbaki labeling of the simple roots $\alpha_i$, as given in \cite[11.4, pg. 58]{Humphreys}. We denote the dominant weights with respect to $\Delta$ by $X(T)^+$, and the fundamental dominant weight corresponding to $\alpha_i$ is denoted by $\omega_i$. We set the usual partial ordering $\preceq$ on $X(T)$, i.e. for $\mu, \lambda \in X(T)$ we have $\mu \preceq \lambda$ if and only if $\lambda = \mu$, or $\lambda - \mu$ is a sum of positive roots. For a dominant weight $\lambda \in X(T)^+$, we denote by $L_G(\lambda)$ the irreducible $G$-module with highest weight $\lambda$, by $V_G(\lambda)$ the Weyl module of highest weight $\lambda$, and by $T_G(\lambda)$ the indecomposable tilting module of highest weight $\lambda$. The longest element in the Weyl group of $G$ is denoted by $w_0$. The character of a $G$-module $V$ is denoted by $\ch V$, and it is the element of $\Z [X(T)]$ defined by $$\ch V = \sum_{\mu \in X(T)} m_V(\mu) \mu,$$ where $m_V(\mu)$ is the dimension of the $\mu$-weight space of $V$. 

Let $F: G \rightarrow G$ be the Frobenius endomorphism induced by the field automorphism $x \mapsto x^p$ of $\fieldsymbol$, see for example \cite[Lemma 76]{SteinbergNotes}. For a $G$-module $V$ corresponding to the representation $\rho: G \rightarrow \GL(V)$, we denote the $G$-module corresponding to the representation $\rho \circ F^k: G \rightarrow \GL(V)$ by $V^{[k]}$. The $G$-module $V^{[k]}$ is called the $k$th Frobenius twist of $V$.

The socle of a $G$-module $V$ is denoted by $\soc V$. If a $G$-module $V$ has a filtration $V = V_1 \supset V_2 \supset \cdots \supset V_{t} \supset V_{t+1} = 0$ with $W_i \cong V_i / V_{i+1}$, we will denote this by $V = W_1 | W_2 | \cdots | W_t$.

Throughout $V$ will always denote a finite-dimensional vector space over $K$. Let $u \in \GL(V)$ be a unipotent linear map. It will often be convenient for us to describe the action of $u$ on a representation in terms of $\fieldsymbol [u]$-modules. Suppose that $u$ has order $q = p^t$. Then there exist exactly $q$ indecomposable $\fieldsymbol [u]$-modules which we will denote by $J_1$, $J_2$, $\ldots$, $J_q$. Here $\dim J_i = i$ and $u$ acts on $J_i$ as a full Jordan block. We use the notation $r \cdot J_n$ for the direct sum $J_n \oplus \cdots \oplus J_n$, where $J_n$ occurs $r$ times.

A bilinear form $b$ on $V$ is \emph{non-degenerate}, if its \emph{radical} $\rad b = \{v \in V: b(v,w) = 0 \text{ for all } w \in V \}$ is zero. For a quadratic form $Q : V \rightarrow \fieldsymbol$ on a vector space $V$, its \emph{polarization} is the bilinear form $b_Q$ defined by $b_Q(v,w) = Q(v+w) - Q(v) - Q(w)$ for all $v, w \in V$.  We say that $Q$ is \emph{non-degenerate}, if its \emph{radical} $\rad Q = \{v \in \rad b_Q : Q(v) = 0 \}$ is zero.

\section{Preliminaries on unipotent elements}

For a unipotent linear map $u \in \GL(V)$ of order $p$, write $r_p(u)$ for the number of Jordan blocks of size $p$ in the Jordan decomposition of $u$. Using the fact that $r_p(u) = \rank (u-1)^{p-1}$, it is easy to prove the following lemma, as observed in \cite[(1), pg. 2585]{SuprunenkoAsymptotic}.

\begin{lemma}\label{jordanblockfiltration}
Let $u$ be a unipotent linear map on the vector space $V$ and suppose that $u$ has order $p$. Suppose that $V$ has a filtration $V = W_1 \supseteq W_2 \supseteq \cdots \supseteq W_t \supseteq W_{t+1} = 0$ of $K [u]$-submodules. Then $r_p(u) \geq \sum_{i = 1}^t r_p(u_{W_i/W_{i+1}})$.
\end{lemma}

We will need the following result to decompose tensor products of unipotent matrices of order $p$.

\begin{lemma}[{\cite[Theorem 1]{Renaud}}]\label{lemmatensordecomp_orderp}
Let $1 \leq m \leq n \leq p$. Then $$J_m \otimes J_n \cong \oplus_{i = 0}^{h-1} J_{n-m+2i+1} \oplus N \cdot J_p$$ where $h = \min \{m, p-n\}$ and $N = \max \{0, m+n-p\}$. In particular, $J_m \otimes J_p = m \cdot J_p$ for all $1 \leq m \leq p$.
\end{lemma}

We state next some results on the distinguished unipotent conjugacy classes in the simple classical groups $\SL(V)$, $\Sp(V)$, and $\SO(V)$. In good characteristic, the distinguished unipotent classes are described by the next three lemmas. For proofs, see for example \cite[Proposition 3.5]{LiebeckSeitzClass}.

\begin{lemma}\label{lemma:dist_SL}
Let $u \in \SL(V)$ be a unipotent element. Then $u$ is a distinguished unipotent element of $\SL(V)$ if and only if $V \downarrow K[u] = J_d$, where $d = \dim V$.
\end{lemma}

\begin{lemma}\label{lemma:dist_odd_SP}
Assume $p \neq 2$. Let $u \in \Sp(V)$ be a unipotent element. Then $u$ is a distinguished unipotent element of $\Sp(V)$ if and only if $$V \downarrow K[u] = J_{d_1} \oplus \cdots \oplus J_{d_t},$$ where $d_i$ are distinct even integers.
\end{lemma}

\begin{lemma}\label{lemma:dist_odd_SO}
Assume $p \neq 2$. Let $u \in \SO(V)$ be a unipotent element. Then $u$ is a distinguished unipotent element of $\SO(V)$ if and only if $$V \downarrow K[u] = J_{d_1} \oplus \cdots \oplus J_{d_t},$$ where $d_i$ are distinct odd integers.
\end{lemma}

In characteristic two, a classification of the unipotent conjugacy classes of $\Sp(V)$ and $\SO(V)$ was given by Hesselink in \cite{Hesselink}. Starting from this result, Liebeck and Seitz have given detailed information about the structure of the centralizers of unipotent elements in \cite{LiebeckSeitzClass}. In particular, they have given a description of the distinguished unipotent classes, which we record in the next lemma.

\begin{lemma}[{\cite[Proposition 6.1, Section 6.8]{LiebeckSeitzClass}}]\label{lemma:dist_even_SP}
Assume $p = 2$. Let $G = \Sp(V)$ or $G = \SO(V)$. Set $Z = V$ if $\dim V$ is even and $Z = V/V^\perp$ if $\dim V$ is odd. Let $u \in G$ be a unipotent element. Then $u$ is a distinguished unipotent element of $G$ if and only if there is an orthogonal decomposition $$Z \downarrow K[u] = J_{d_1} \perp J_{d_2} \perp \cdots \perp J_{d_t},$$ where $d_i$ is even for all $1 \leq i \leq t$, and each $J_{d_i}$ occurs with multiplicity at most two.
\end{lemma}


\begin{lemma}\label{lemma:dist_even_SP_noreps}
Assume $p = 2$. Let $G = \Sp(V)$ or $G = \SO(V)$, where $\dim V$ is even. Let $u \in G$ be a unipotent element with Jordan form $$V \downarrow K[u] = J_{d_1} \oplus \cdots \oplus J_{d_t},$$ where $d_i$ are distinct even integers. Then there exists an orthogonal decomposition $V \downarrow K[u] = J_{d_1} \perp J_{d_2} \perp \cdots \perp J_{d_t}$, and $u$ is a distinguished unipotent element of $G$.
\end{lemma}

\begin{proof}This is an easy consequence of the \emph{distinguished normal form} for $u$ established in \cite[Chapter 6]{LiebeckSeitzClass}. In the orthogonal decomposition given in \cite[Lemma 6.2]{LiebeckSeitzClass}, there cannot be any summands of the form $W(m)$, since $u$ acts on $W(m)$ with two Jordan blocks of size $m$. Thus the claim follows from \cite[Proposition 6.1]{LiebeckSeitzClass}.\end{proof}

\begin{lemma}\label{lemma:typeBtoDstabilizer}
Assume $p = 2$. Let $G = \SO(V, Q)$, where $\dim V$ is even and $Q$ is a non-degenerate quadratic form on $V$ with polarization $\beta$. Fix a vector $v \in V$ such that $Q(v) \neq 0$. Suppose that $u \in \operatorname{Stab}_G(v)$ is a unipotent element such that with respect to the alternating bilinear form induced on $\langle v \rangle^\perp / \langle v \rangle$ by $\beta$, we have an orthogonal decomposition $$\langle v \rangle^\perp / \langle v \rangle \downarrow K[u] = J_{d_1} \perp \cdots \perp J_{d_t}.$$ Then $$V \downarrow K[u] = \begin{cases} (J_1 \oplus J_1) \perp J_{d_1} \perp \cdots \perp J_{d_t}, & \text{ if } t \text{ is even,} \\ J_2 \perp J_{d_1} \perp \cdots \perp J_{d_t}, & \text{ if } t \text{ is odd.} \end{cases}$$
\end{lemma}

\begin{proof}It follows from \cite[6.8]{LiebeckSeitzClass} that there exists an orthogonal decomposition $\langle v \rangle^\perp \downarrow K[u] = W \perp \langle v \rangle$ such that $W \downarrow K[u] = J_{d_1} \perp \cdots \perp J_{d_t}$. Now $W$ is a non-degenerate subspace, so we have an orthogonal decomposition $V \downarrow K[u] = W \perp Z$, where $\dim Z = 2$. It is obvious that $Z \downarrow K[u] = J_1 \oplus J_1$ or $Z \downarrow K[u] = J_2$. On the other hand, we have $u \in \SO(V)$, so the number of Jordan blocks of $u$ must be even \cite[Proposition 6.22 (i)]{LiebeckSeitzClass}. Consequently $Z \downarrow K[u] = J_1 \oplus J_1$ if $t$ is even, and $Z \downarrow K[u] = J_2$ if $t$ is odd, as claimed.\end{proof}



%% file: PROJECT_quasiminuscule.tex
%

\begin{table}[!htbp]
\centering
\begin{tabular}{| l | l | l  l |}
\hline
Type               & $\lambda$ & Structure of $V_G(\lambda)$ & Conditions \\
\hline
$A_l$, $l \geq 1$  & $\omega_1 + \omega_l$    & $L_G(\lambda)$ & $p \nmid l+1$ \\
                   &                          & $L_G(\lambda) | L_G(0)$ & $p \mid l+1$ \\
& & & \\									
$B_l$, $l \geq 2$  & $\omega_1$ & $L_G(\lambda)$ & $p \neq 2$ \\
                   &            & $L_G(\lambda) | L_G(0)$ & $p = 2$ \\
& & & \\									
$C_l$, $l \geq 2$  & $\omega_2$ & $L_G(\lambda)$ & $p \nmid l$ \\
                   &            & $L_G(\lambda) | L_G(0)$ & $p \mid l$ \\
& & & \\							
$D_l$, $l \geq 4$  & $\omega_2$ & $L_G(\lambda)$ & $p \neq 2$ \\
                   &            & $L_G(\lambda) | L_G(0)$ & $p = 2$, $l$ odd \\
									 &						& $L_G(\lambda) | L_G(0)^2$ & $p = 2$, $l$ even \\
& & & \\								
$G_2$ & $\omega_1$ & $L_G(\lambda)$ & $p \neq 2$ \\
      &            & $L_G(\lambda) | L_G(0)$ & $p = 2$ \\
& & & \\		
$F_4$ & $\omega_4$ & $L_G(\lambda)$ & $p \neq 3$ \\
      &            & $L_G(\lambda) | L_G(0)$ & $p = 3$ \\
& & & \\			
$E_6$ & $\omega_2$ & $L_G(\lambda)$ & $p \neq 3$ \\
      &            & $L_G(\lambda) | L_G(0)$ & $p = 3$ \\
& & & \\			
$E_7$ & $\omega_1$ & $L_G(\lambda)$ & $p \neq 2$ \\
      &            & $L_G(\lambda) | L_G(0)$ & $p = 2$ \\
& & & \\			
$E_8$ & $\omega_8$ & $L_G(\lambda)$ & none \\
\hline
\end{tabular}
\caption{Quasi-minuscule weights $\lambda \in X(T)^+$.}\label{table:quasiminuscule}
\end{table}

We say that a non-zero dominant weight $\lambda \in X(T)^+$ is \emph{quasi-minuscule}, if $\lambda \succ 0$ and the only weights subdominant to $\lambda$ are $0$ and $\lambda$ itself. This is equivalent to saying that $\lambda \succ 0$ and all non-zero weights occurring in $V_G(\lambda)$ are conjugate under the action of the Weyl group of $G$. It is well known that there exists a unique quasi-minuscule weight $\lambda \in X(T)^+$, and it is equal to the highest short root of $G$. We give $\lambda$ explicitly in Table \ref{table:quasiminuscule}. In Table \ref{table:quasiminuscule} we have also given the structure of the corresponding Weyl module $V_G(\lambda)$, see for example \cite[Theorem 5.1]{Lubeck}.

In this section, we give some results on the structure of indecomposable tilting modules $T_G(\lambda)$ with quasi-minuscule highest weight $\lambda \in X(T)^+$. We know that $T_G(\lambda) = V_G(\lambda)$ if $V_G(\lambda)$ is irreducible. If $V_G(\lambda)$ is not irreducible and $G$ is not of type $D_l$, then $V_G(\lambda) = L_G(\lambda) | L_G(0)$ and in this case we will see that $T_G(\lambda) = L_G(0) | L_G(\lambda) | L_G(0)$ (Lemma \ref{lemma:forms_quasiminuscule} (i)). When $V_G(\lambda) = L_G(\lambda) | L_G(0)$, we will also establish results on the existence of non-degenerate $G$-invariant forms on $T_G(\lambda)$ (Lemma \ref{lemma:forms_quasiminuscule} (ii)-(iii)). In Section \ref{section:nonGcr_orderp}, we will apply these results to present examples of non-completely reducible subgroups of classical groups containing distinguished unipotent elements.


We begin with the following lemma. From the proof it is clear that nothing specific to algebraic groups is needed, and indeed the result is true in a more general setting.

\begin{lemma}\label{lemma:char2_symmetric_to_alternating}
Assume $p = 2$. Let $V$ be an indecomposable $G$-module. If $V$ admits a non-degenerate $G$-invariant symmetric bilinear form, then $V$ admits a non-degenerate $G$-invariant alternating bilinear form.
\end{lemma}

\begin{proof}
Let $\beta$ be a non-degenerate $G$-invariant symmetric bilinear form on $V$. Since $p = 2$, the map $f: V \rightarrow K$ defined by $f(v) = \sqrt{\beta(v,v)}$ is a morphism of $G$-modules. If $f = 0$, then $\beta$ is alternating and we are done. 

Suppose then that $f \neq 0$. The bilinear form $\beta$ is non-degenerate, so there exists a non-zero $z \in V$ such that $f(v) = \beta(v,z)$ for all $v \in V$. Since $f$ is a $G$-morphism, the vector $z$ is a $G$-fixed point in $V$. Hence the subspace $\langle z \rangle$ must be totally isotropic with respect to $\beta$, since $V$ is indecomposable. That is, we have $\beta(z,z) = f(z) = 0$.

Consider the map $N: V \rightarrow V$ defined by $N(v) = \beta(v,z)z = f(v)z$ for all $v \in V$. Then $N \in \End_G(V)$ since $z$ is fixed by $G$, and $N^2 = 0$ since $f(z) = 0$. Hence the map $\psi = \operatorname{id}_V + N$ is an isomorphism of $G$-modules. This gives a non-degenerate $G$-invariant symmetric bilinear form $\gamma$ on $V$ via $$\gamma(v,v') = \beta(\psi(v), v') = \beta(v,v') + f(v)f(v')$$ for all $v,v' \in V$. The bilinear form $\gamma$ is alternating, since $\gamma(v,v) = \beta(v,v) + f(v)^2 = 0$ for all $v \in V$.\end{proof}

\begin{lemma}\label{lemma:forms_quasiminuscule}
Let $\lambda \in X(T)^+$ be quasi-minuscule. Suppose that $V_G(\lambda) = L_G(\lambda) | L_G(0)$. Then:

\begin{enumerate}[\normalfont (i)]
\item The indecomposable tilting module $T_G(\lambda)$ is uniserial, and $T_G(\lambda) = L_G(0) | L_G(\lambda) | L_G(0)$. 
\item If $p \neq 2$, then $T_G(\lambda)$ admits a non-degenerate $G$-invariant symmetric bilinear form.
\item If $p = 2$, then $T_G(\lambda)$ admits a non-degenerate $G$-invariant alternating bilinear form, unique up to a scalar multiple.
\end{enumerate}
\end{lemma}

\begin{proof}
For (i), note that from the short exact sequence $$0 \rightarrow L_G(0) \rightarrow V_G(\lambda) \rightarrow L_G(\lambda) \rightarrow 0$$ the long exact sequence in cohomology gives $H^1(G, V_G(\lambda)) \cong H^1(G, L_G(\lambda))$, since $H^1(G, L_G(0)) = 0 = H^2(G, L_G(0))$ by \cite[Corollary II.4.11]{JantzenBook}. On the other hand, by the assumption $V_G(\lambda) = L_G(\lambda) | L_G(0)$ and \cite[Proposition II.2.14]{JantzenBook} we have $H^1(G, L_G(\lambda)) \cong K$. Therefore $H^1(G, V_G(\lambda)) \cong K$, and so up to isomorphism there exists a unique nonsplit extension $$0 \rightarrow V_G(\lambda) \rightarrow V \rightarrow L_G(0) \rightarrow 0.$$ The fact that $V$ is a nonsplit extension of $L_G(0)$ by $V_G(\lambda)$ implies that $\operatorname{soc} V \cong \operatorname{soc} V_G(\lambda) \cong K$. Consequently $V / \soc V$ is a nonsplit extension of $L_G(0)$ by $L_G(\lambda)$, so it follows that $V$ is uniserial and $V / \soc V \cong V_G(-w_0\lambda)^*$. Hence $V$ has a filtration by dual Weyl modules, and by construction $V$ has a filtration by Weyl modules. Thus $V$ is a tilting module, so $V \cong T_G(\lambda)$ since $V$ is indecomposable with highest weight $\lambda$. This establishes (i).

For claims (ii) and (iii), note first that the fact that $\lambda$ is quasi-minuscule implies that $-w_0 \lambda = \lambda$, where $w_0$ is the longest element of the Weyl group of $G$. Thus $T_G(\lambda)$ is self-dual, see \cite[Remark II.E.6]{JantzenBook}. 

We suppose first that $p \neq 2$ and consider claim (ii). Since $T_G(\lambda)$ is indecomposable, its endomorphism ring is a local ring. It follows then from a general result in representation theory \cite[Lemma 2.1]{Quebbemann} (or \cite[Satz 2.11 (a)]{WillemsThesis}) that there exists a non-degenerate $G$-invariant bilinear form $\beta$ on $T_G(\lambda)$ such that $\beta$ is symmetric or alternating. Since $T_G(\lambda)$ is uniserial and $T_G(\lambda) = L_G(0) | L_G(\lambda) | L_G(0)$, the form $\beta$ induces a non-degenerate $G$-invariant bilinear form on the subquotient of $T_G(\lambda)$ isomorphic to $L_G(\lambda)$. Since $\lambda \succ 0$, it follows from \cite[Lemma 79, pg. 226-227]{SteinbergNotes} that the form induced on the subquotient is symmetric. Because $p \neq 2$, we conclude that $\beta$ must also be symmetric, which establishes (ii).

For claim (iii), suppose that $p = 2$ and set $V = T_G(\lambda)$. We first show that there exists a non-degenerate $G$-invariant alternating bilinear form on $V$. Fix some isomorphism $f: V \rightarrow V^*$ of $G$-modules. This gives a non-degenerate $G$-invariant bilinear form $\beta$ on $V$ via $\beta(v,v') = f(v)(v')$ for all $v,v' \in V$. Consider the bilinear form $\gamma$ on $V$ defined by $\gamma(v,v') = \beta(v,v') + \beta(v',v)$ for all $v,v' \in V$. If $\gamma = 0$, then $\beta$ is symmetric and Lemma \ref{lemma:char2_symmetric_to_alternating} gives a non-degenerate $G$-invariant alternating bilinear form on $V$. Suppose then that $\gamma \neq 0$. Then $\gamma$ is a non-zero $G$-invariant alternating bilinear form on $V$, and we will show that $\gamma$ must be non-degenerate. Note that $\gamma$ induces a non-degenerate $G$-invariant alternating bilinear form on $V / \operatorname{rad} \gamma$. Since $V$ is uniserial, and since $L_G(0)$ and $V_G(\lambda)^* = L_G(0) | L_G(\lambda)$ do not admit non-degenerate $G$-invariant alternating bilinear forms, the only possibility is that $\operatorname{rad} \gamma = 0$. In other words, the bilinear form $\gamma$ is non-degenerate.

Next we show the uniqueness statement of (iii). Let $\beta$ and $\gamma$ be non-degenerate $G$-invariant alternating bilinear forms on $V$. It is easy to see (for example \cite[Satz 2.3]{WillemsThesis}) that there exists a unique $G$-isomorphism $\varphi \in \operatorname{End}_G(V)$ such that $\gamma(v,v') = \beta(\varphi(v), v')$ for all $v, v' \in V$. 

We now describe $\operatorname{End}_G(V)$. We have established in the proof of (i) that $\Ext_G^1(L_G(0), L_G(\lambda)) \cong K$ and that $T_G(\lambda) = L_G(0) | L_G(\lambda) | L_G(0)$ is uniserial, so by \cite[Example 4.5]{GaribaldiNakano} we have $\dim \operatorname{End}_G(V) = 2$. Let $z \in V$ be a vector spanning the unique $1$-dimensional submodule of $V$. As a basis for the vector space $\operatorname{End}_G(V)$, we can take the identity map $1: V \rightarrow V$ and the nilpotent map $N: V \rightarrow V$ defined by $N(v) = \beta(v,z)z$ for all $v \in V$. 

Write $\varphi = \lambda \cdot 1 + \mu \cdot N$. Since $\gamma$ and $\beta$ are alternating, we have $\gamma(v,v) = \mu \cdot \beta(Nv, v) = \mu \cdot \beta(v,z)^2$ for all $v \in V$. Now since $\beta$ is non-degenerate, we can choose a $v \in V$ with $\beta(v,z) \neq 0$, and consequently $\mu = 0$. Therefore $\gamma = \lambda \cdot \beta$, which completes the proof of (iii) and the lemma.\end{proof}

\begin{remark}
Lemma \ref{lemma:forms_quasiminuscule} (ii) is true much more generally. Suppose that $p \neq 2$. If $\lambda = -w_0 \lambda$, then $T_G(\lambda)$ is self-dual and one can show that $T_G(\lambda)$ admits a non-degenerate $G$-invariant bilinear form $\beta$ which is alternating or symmetric. Furthermore, arguing similarly to the proof of Lemma \ref{lemma:forms_quasiminuscule} (ii), one can show that $\beta$ must be of the same type as the correponding form on $L_G(\lambda)$.
\end{remark}

\begin{remark}\label{remark:quasi_typeD_omitted}
In Table \ref{table:quasiminuscule} we have given the cases where Lemma \ref{lemma:forms_quasiminuscule} applies. For $T_G(\lambda)$ with quasi-minuscule highest weight $\lambda \in X(T)^+$, we have left untreated the case where $G$ is simple of type $D_l$ when $p = 2$ and $l$ is even. In this case we have $V_G(\lambda) = L_G(\lambda) | L_G(0)^2$. We omit the proof, but proceeding similarly to the proof of Lemma \ref{lemma:forms_quasiminuscule}, one can show that $T_G(\lambda) = L_G(0)^2 | L_G(\lambda) | L_G(0)^2$ and that $T_G(\lambda)$ admits a non-degenerate $G$-invariant alternating bilinear form. In this case such a form is not unique up to a scalar multiple, but it is unique up to a $G$-equivariant isometry.
\end{remark}

%% file: PROJECT_SL2.tex
%

In this section, let $G$ be the algebraic group $\SL_2(K)$ with natural module $E$. Fix also a nonidentity unipotent element $u \in G$. Throughout we will identify the weights of a maximal torus of $G$ with $\Z$, and the dominant weights with $\Z_{\geq 0}$. With this identification, for weights $\lambda, \mu \in \Z$ we have $\mu \preceq \lambda$ if and only if $\lambda - \mu$ is a non-negative integer such that $\lambda \equiv \mu \mod{2}$.

The main purpose of this section is to classify indecomposable $G$-modules $V$ where $u$ acts on $V$ with at most one Jordan block of size $p$. One consequence of this is a criterion for a representation of $G$ to be semisimple. Specifically, we prove that a self-dual $G$-module $V$ must be semisimple if $u$ acts on $V$ with at most one Jordan block of size $p$ (Proposition \ref{SL2_semisimplicity}).

We begin by stating two well known results in the representation theory of $G$, which we will use throughout this section. Below for a positive integer $n$, we use the notation $\nu_p(n)$ for the largest integer $k \geq 0$ such that $p^k$ divides $n$.

\begin{lause}[{Cline, \cite[Corollary 3.9]{Landrock}}]\label{sl2split}
Let $\lambda = \sum_{i \geq 0} \lambda_i p^i$ and $\mu = \sum_{i \geq 0} \mu_i p^i$ be weights of $G$, where $0 \leq \lambda_i, \mu_i < p$. Then $\Ext_G^1(L_G(\lambda), L_G(\mu)) \neq 0$ precisely when there exists a $k \geq \nu_p(\lambda + 1)$ such that $\mu_i = \lambda_i$ for all $i \neq k,k+1$ and $\mu_k = p-2-\lambda_k$, $\mu_{k+1} = \lambda_{k+1} \pm 1$. In the case where $\Ext_G^1(L_G(\lambda), L_G(\mu)) \neq 0$, we have $\Ext_G^1(L_G(\lambda), L_G(\mu)) \cong K$.
\end{lause}


\begin{lause}[{\cite[Lemma 2.3]{SeitzSaturation}, \cite[Example 2, pg. 47]{DonkinTilting}}]\label{thm:tiltingSL2}
Let $c \in \Z_{\geq 0}$ be a dominant weight. Then:

\begin{enumerate}[\normalfont (i)]
\item If $0 \leq c \leq p-1$, the indecomposable tilting module $T_G(c)$ is irreducible, so $T_G(c) = L_G(c)$.
\item If $p \leq c \leq 2p-2$, the indecomposable tilting module $T_G(c)$ is uniserial of dimension $2p$, and $T_G(c) = L_G(2p-2-c) | L_G(c) | L_G(2p-2-c)$.
\item (Donkin) If $c > 2p-2$, then $T_G(c) \cong T_G(p-1+r) \otimes T_G(s)^{[1]}$, where $s \geq 1$ and $0 \leq r \leq p-1$ are such that $c = sp + (p-1+r)$.
\end{enumerate}
\end{lause}


We next describe the Jordan block sizes of $u$ acting on Weyl modules and tilting modules of $G$.

\begin{lemma}\label{symmetricpowerlemma}
Let $m \in \Z_{\geq 0}$ and write $m = qp + r$, where $q \geq 0$ and $0 \leq r < p$. Then $V_G(m) \downarrow K[u] = q \cdot J_p \oplus J_{r+1}$.
\end{lemma}

\begin{proof}
The Weyl module $V_G(m)$ is isomorphic to the dual of the symmetric power $S^m(E)$, see for example \cite[II.2.16]{JantzenBook}. Therefore it suffices to compute the Jordan block sizes of $u$ acting on the symmetric power $S^m(E)$.

Fix a basis $x, y$ of $E$ such that $ux = x$ and $uy = x+y$. Consider the basis $x^m, x^{m-1}y, \ldots, y^m$ induced on $S^m(E)$. Then $$u(x^{m-k}y^k) = x^{m-k}(x+y)^k = \sum_{i = 0}^k \binom{k}{i} x^{m-i}y^i,$$ so with respect to this basis, the matrix of $u$ acting on $S^m(E)$ is the upper triangular Pascal matrix $P = ( \binom{i-1}{j-1} )_{0 \leq i,j \leq m}$. Here we define $\binom{i}{j} = 0$ if $i < j$. The Jordan form of the transpose of $P$ is computed for example in \cite{Callan}, and from this result we find that the transpose of $P$ has $q$ Jordan blocks of size $p$ and one Jordan block of size $r+1$. Since a matrix is similar to its transpose, the lemma follows.\end{proof}

\begin{lemma}\label{lemma:tiltingSL2_jordanblocks}
Let $c \geq p-1$. Then $T_G(c) \downarrow K [u] = N \cdot J_p$, where $N = \dim T_G(c) / p$. In particular, $T_G(c) \downarrow K [u] = J_p \oplus J_p$ if $p \leq c \leq 2p-2$.
\end{lemma}

\begin{proof}We argue similarly to \cite[Proposition 5]{McNinchAdjoint}. Without loss of generality, we can assume that $u$ is contained in the finite subgroup $G(p) := \SL_2(\mathbb{F}_p)$ of $G$. We note first that for all $p-1 \leq c \leq 2p-2$, the restriction of $T_G(c)$ to $G(p)$ is projective. For $c = p-1$ this follows since $T_G(c)$ is the Steinberg module, and for $p \leq c \leq 2p-2$ this follows from \cite[Lemma 2.3 (c)]{SeitzSaturation}. 

A tensor product with a projective module remains projective \cite[Lemma 4, pg. 47]{Alperin}, so we conclude from Theorem \ref{thm:tiltingSL2} (iii) that the restriction of $T_G(c)$ to $G(p)$ is projective for all $c \geq p-1$. Then the restriction of $T_G(c)$ to the Sylow $p$-subgroup $H = \langle u \rangle$ of $G(p)$ is also projective \cite[Theorem 6, pg. 33]{Alperin}. Since $J_p$ is the only projective indecomposable $H$-module, the claim follows.\end{proof}

%

%

\begin{prop}\label{nonsplitaction}
Let $V$ be a $G$-module which is a nonsplit extension of $L_G(\lambda)$ by $L_G(\mu)$. Then $u$ acts on $V$ with at least one Jordan block of size $p$. If $u$ acts on $V$ with precisely one Jordan block of size $p$, then there exists $p \leq c \leq 2p-2$ and $l \geq 0$ such that one of the following holds:

\begin{enumerate}[\normalfont (i)]
\item $\lambda = cp^l$, $\mu = (2p-2-c)p^l$ and $V \cong V_G(c)^{[l]}$.
\item $\lambda = (2p-2-c)p^l$, $\mu = cp^l$ and $V \cong (V_G(c)^*)^{[l]}$.
\end{enumerate}

\noindent Moreover, if $V$, $\lambda$, $\mu$, and $c$ are as in case (i) or case (ii), then $V$ is a nonsplit extension of $L_G(\lambda)$ by $L_G(\mu)$, and $V \downarrow K[u] = J_p \oplus J_{c-p+1}$.


\end{prop}

\begin{proof}

We begin by considering the claims about $V_G(c)^{[l]}$ and $(V_G(c)^*)^{[l]}$, where $p \leq c \leq 2p-2$ and $l \geq 0$. It is well known (and easily seen by considering the weights in $V_G(c)$) that $V_G(c)$ is a nonsplit extension $$0 \rightarrow L_G(2p-2-c) \rightarrow V_G(c) \rightarrow L_G(c) \rightarrow 0.$$ Since the irreducible representations of $G$ are self-dual, we see that $V_G(c)^*$ is a nonsplit extension $$0 \rightarrow L_G(c) \rightarrow V_G(c)^* \rightarrow L_G(2p-2-c) \rightarrow 0.$$ Let $\lambda = cp^l$ and $\mu = (2p-2-c)p^l$, where $l \geq 0$. By taking a Frobenius twist, we see that for all $l \geq 0$ the Weyl module $V_G(c)^{[l]}$ is a non-split extension of $L_G(\lambda)$ by $L_G(\mu)$, and its dual $(V_G(c)^*)^{[l]}$ is a non-split extension of $L_G(\mu)$ by $L_G(\lambda)$. Furthermore, by Lemma \ref{symmetricpowerlemma} the unipotent element $u \in G$ acts on $V_G(c)$ with Jordan form $J_p \oplus J_{c-p+1}$. The Jordan block sizes are not changed by taking a dual or a Frobenius twist, so for all $l \geq 0$ the element $u$ acts on both $V_G(c)^{[l]}$ and $(V_G(c)^*)^{[l]}$ with Jordan form $J_p \oplus J_{c-p+1}$. This completes the proof that the properties of $V_G(c)^{[l]}$ and $(V_G(c)^*)^{[l]}$ are as claimed.


Now let $\lambda, \mu \in \Z_{\geq 0}$ be arbitrary weights, and let $V$ be a nonsplit extension $$0 \rightarrow L_G(\mu) \rightarrow V \rightarrow L_G(\lambda) \rightarrow 0.$$ Assume that $u$ acts on $V$ with at most one Jordan block of size $p$. We note first that to prove the proposition, it will be enough to show that there exist $p \leq c \leq 2p-2$ and $l \geq 0$ such that $\lambda$ and $\mu$ are as in case (i) or (ii) of the claim. Indeed, if this holds, then since $\Ext_G^1(L_G(\lambda), L_G(\mu)) \cong K$ (Theorem \ref{sl2split}), we must have $V \cong V_G(c)^{[l]}$ or $V \cong (V_G(c)^*)^{[l]}$. Furthermore, as seen in the first paragraph, then $u$ acts on $V$ with Jordan form $J_p \oplus J_{c-p+1}$, in particular with exactly one Jordan block of size $p$.



We proceed to show that $\lambda$ and $\mu$ are as claimed, which will prove the proposition. Write $\lambda = \sum_{i \geq 0} \lambda_i p^i$, where $0 \leq \lambda_i \leq p-1$ for all $i$. By the Steinberg tensor product theorem, we have

$$L_G(\lambda) \cong \bigotimes_{i \geq 0} L_G(\lambda_i)^{[i]}.$$

Consider first the case where $\lambda_l = p-1$ for some $l$. Then $u$ acts on $L_G(\lambda_l)^{[l]}$ with a single Jordan block of size $p$. Now the tensor product of a Jordan block of size $p$ with any Jordan block of size $c \leq p$ consists of $c$ Jordan blocks of size $p$ (Lemma \ref{lemmatensordecomp_orderp}), so it follows that $L_G(\lambda) \downarrow K [u] = N \cdot J_p$, where $N = \dim L_G(\lambda) / p$. Since $u$ acts on $V$ with at most one Jordan block of size $p$, by Lemma \ref{jordanblockfiltration} the action of $u$ on $L_G(\lambda)$ can have at most one Jordan block of size $p$. It follows then that $N = 1$, so $\lambda = (p-1)p^l$. Since $V$ is a nonsplit extension, by Theorem \ref{sl2split} the weight $\mu$ must be $(p-2)p^{l-1} + (p-2)p^l$ or $(p-1)p^l + (p-2)p^k + p^{k+1}$ for some $k \neq l,l-1$. By Lemma \ref{jordanblockfiltration} the action of $u$ on $L_G(\mu)$ has no Jordan blocks of size $p$. With Lemma \ref{lemmatensordecomp_orderp}, one finds that this happens only if $p = 2$ and $\mu = (p-2)p^{l-1} + (p-2)p^l = 0$. Then $\lambda$ and $\mu$ are as in case (i) of the claim.

Thus we can assume that $0 \leq \lambda_i \leq p-2$ for all $i$. Write $\mu = \sum_{i \geq 0} \mu_i p^i$, where $0 \leq \mu_i \leq p-1$ for all $i$. According to Theorem \ref{sl2split}, there exists an $l$ such that we have $\mu_i = \lambda_i$ for $i \neq l,l+1$, and $\mu_l = p-2-\lambda_l$, $\mu_{l+1} = \lambda_{l+1} \pm 1$. We can write $\lambda = \zeta + \lambda_lp^l + \lambda_{l+1}p^{l+1}$ and $\mu = \zeta + \mu_lp^l + \mu_{l+1}p^{l+1}$. By the Steinberg tensor product theorem, we have
\begin{align*}
L_G(\lambda) &\cong L_G(\zeta) \otimes L_G(\lambda_l + \lambda_{l+1}p)^{[l]} \\
L_G(\mu) &\cong L_G(\zeta) \otimes L_G(\mu_l + \mu_{l+1}p)^{[l]}
\end{align*}

Note that here $p$ does not divide the dimension of $L_G(\zeta)$, because we are assuming that $\lambda_i < p-1$ for all $i$. Furthermore, by Theorem \ref{sl2split} there exists a $G$-module $W$ which is a nonsplit extension $$0 \rightarrow L_G(\lambda_l + \lambda_{l+1}p) \rightarrow W \rightarrow L_G(\mu_l + \mu_{l+1}p) \rightarrow 0.$$ Therefore by \cite[Theorem 2.4]{SerreTensor}, the module $L_G(\zeta) \otimes W^{[l]}$ is a nonsplit extension of $L_G(\lambda)$ by $L_G(\mu)$. Hence $L_G(\zeta) \otimes W^{[l]} \cong V$, because a nonsplit extension of $L_G(\lambda)$ by $L_G(\mu)$ is unique by Theorem \ref{sl2split}. 

We will first treat the case where $\zeta = 0$. Here $V \cong W^{[l]}$, so without loss of generality we may assume that $l = 0$. Write $\lambda = c + dp$ and $\mu = (p-2-c) + (d \pm 1)p$, where $0 \leq c, d \leq p-2$ and $d > 0$ if $\mu = (p-2-c) + (d - 1)p$.

Suppose that $\lambda = c + dp$ and $\mu = (p-2-c) + (d-1)p$. Here $\dim V = \dim L_G(\lambda) + \dim L_G(\mu) = \lambda + 1 = \dim V_G(\lambda)$. Now $V$ is a nonsplit extension of $L_G(\lambda)$ by $L_G(\mu)$ and $\lambda \succ \mu$, so by \cite[Lemma II.2.13 (b)]{JantzenBook} we must have $V \cong V_G(\lambda)$. By Lemma \ref{symmetricpowerlemma}, the action of $u$ on $V$ has Jordan form $d \cdot J_p \oplus J_{c+1}$. Therefore $u$ acts on $V$ with at most one Jordan block of size $p$ if and only if $d = 1$, that is, when $\lambda = c + p$ and $\mu = (p-2-c)$. Then $\lambda$ and $\mu$ are as in case (i) of the claim.

Next consider $\lambda = c + dp$ and $\mu = (p-2-c) + (d+1)p$. In this case $\dim V = \dim L_G(\lambda) + \dim L_G(\mu) = \mu + 1 = \dim V_G(\mu)$. Because $V^*$ is a nonsplit extension of $L_G(\mu)$ by $L_G(\lambda)$ and $\mu \succ \lambda$, by \cite[Lemma II.2.13 (b)]{JantzenBook} we have $V^* \cong V_G(\mu)$. By Lemma \ref{symmetricpowerlemma}, the action of $u$ on $V$ has Jordan form $(d+1) \cdot J_p \oplus J_{p-(c+1)}$. In this case $u$ acts with at most one Jordan block of size $p$ if and only if $d = 0$, that is, when $\lambda = c$ and $\mu = (p-2-c) + p$. Then $\lambda$ and $\mu$ are as in case (ii) of the claim, and this completes the proof of the proposition in the case where $\zeta = 0$.

It remains to consider the possibility that $\zeta \neq 0$. Since $W$ is a nonsplit extension of two irreducibles, it follows from the $\zeta = 0$ case that $u$ acts on $W$ with at least one Jordan block of size $p$. Now the tensor product of a Jordan block of size $p$ with any Jordan block of size $c \leq p$ consists of $c$ Jordan blocks of size $p$ by Lemma \ref{lemmatensordecomp_orderp}. Therefore if $\zeta \neq 0$, then $u$ acts on $V$ with more than one Jordan block of size $p$, contradiction.\end{proof}

\begin{lemma}\label{someextlemma}
Let $p \leq c \leq 2p-2$. Then for all $l \geq 0$:
\begin{enumerate}[\normalfont (i)]
\item $\Ext_G^1(V_G(c)^{[l]}, L_G(c)^{[l]}) = 0$.
\item $\Ext_G^1(V_G(c)^{[l]}, L_G(2p-2-c)^{[l]}) = 0$.
\item $\Ext_G^1(L_G(c)^{[l]}, V_G(c)^{[l]}) = 0$.
\item $\Ext_G^1(L_G(2p-2-c)^{[l]}, V_G(c)^{[l]}) \cong K$.
\end{enumerate}

Furthermore, every nonsplit extension of $L_G(2p-2-c)^{[l]}$ by $V_G(c)^{[l]}$ is isomorphic to $T_G(c)^{[l]}$.
\end{lemma}

\begin{proof}Set $c' = 2p-2-c$. We note first that the last claim of the lemma follows from (iv), once we show that $T_G(c)^{[l]}$ is a nonsplit extension of $L_G(c')^{[l]}$ by $V_G(c)^{[l]}$. To this end, by \cite[Lemma 2.3 (b)]{SeitzSaturation} the tilting module $T_G(c)$ is a nonsplit extension of $L_G(c')$ by $V_G(c)$. The extension stays nonsplit after taking a Frobenius twist, so $T_G(c)^{[l]}$ is a nonsplit extension of $L_G(c')^{[l]}$ by $V_G(c)^{[l]}$.

For claims (i) - (iv), we prove them first in the case where $l = 0$. In this case, claims (i) and (ii) follow from the fact $\Ext_G^1(V_G(\lambda), L_G(\mu)) = 0$ for any $\mu \preceq \lambda$ \cite[II.2.14]{JantzenBook}. For (iii) and (iv), recall that there is an exact sequence $$0 \rightarrow L_G(c') \rightarrow V_G(c) \rightarrow L_G(c) \rightarrow 0.$$ Applying the functor $\Hom_G(L_G(d), -)$ gives a long exact sequence \begin{align*}0 &\rightarrow \Hom_G(L_G(d), L_G(c')) \rightarrow \Hom_G (L_G(d), V_G(c)) \rightarrow \Hom_G(L_G(d), L_G(c)) \\ 
& \rightarrow \Ext_G^1(L_G(d), L_G(c')) \rightarrow \Ext_G^1(L_G(d), V_G(c)) \rightarrow \Ext_G^1(L_G(d), L_G(c)) \\ &\rightarrow \cdots \end{align*} Considering this long exact sequence with $d = c$, we get an exact sequence $$0 \rightarrow K  \rightarrow K  \rightarrow \Ext_G^1(L_G(c), V_G(c)) \rightarrow 0$$ since $\Ext_G^1(L_G(c), L_G(c')) \cong K$ and $\Ext_G^1(L_G(c), L_G(c)) = 0$ by Theorem \ref{sl2split} and \cite[II.2.12]{JantzenBook}, respectively. This proves (iii).

With $d = 2p-2-c$, we get an exact sequence $$0 \rightarrow \Ext_G^1(L_G(c'), V_G(c)) \rightarrow K$$ and so $\dim \Ext_G^1(L_G(c'), V_G(c)) \leq 1$. Thus to prove (iv), it will be enough to show that there exists some nonsplit extension of $L_G(c')$ by $V_G(c)$. For this, we have already noted in the beginning of the proof that $T_G(c)$ is such an extension.

We consider then claims (i) - (iv) for $l > 0$. If $p \neq 2$, then the claims follow from the case $l = 0$, since $\Ext_G^1(X^{[l]}, Y^{[l]}) \cong \Ext_G^1(X, Y)$ for all $G$-modules $X$ and $Y$ \cite[II.10.17]{JantzenBook}. Suppose then that $p = 2$. Note that then we must have $c = 2$, and $c' = 0$. By \cite[Proposition II.10.16, Remark 12.2]{JantzenBook}, for any $G$-modules $X$ and $Y$ there exists a short exact sequence \begin{align*}
0 \rightarrow \Ext_G^1 & (X, Y) \rightarrow \Ext_G^1(X^{[l]}, Y^{[l]}) \\
&\rightarrow \Hom_G(X, L_G(1) \otimes Y) \rightarrow 0.
\end{align*}

Thus claims (i) - (iv) will follow from the case $l = 0$ once we show that $\Hom_G(X, L_G(1) \otimes Y) = 0$ in the following cases:

\begin{enumerate}[(i)$'$]
\item $X = V_G(2)$ and $Y = L_G(2)$,
\item $X = V_G(2)$ and $Y = L_G(0)$,
\item $X = L_G(2)$ and $Y = V_G(2)$,
\item $X = L_G(0)$ and $Y = V_G(2)$.
\end{enumerate}

In all cases (i)$'$ - (iv)$'$, it is straightforward to see that $X$ and $L_G(1) \otimes Y$ have no composition factors in common. Thus $\Hom_G(X, L_G(1) \otimes Y) = 0$, as claimed.\end{proof} 

%
%


We are now ready to prove the main results of this section.

\begin{prop}\label{fewblocksofsizep}
Let $V$ be an indecomposable $G$-module. Then one of the following holds:
\begin{enumerate}[\normalfont (i)]
\item $u$ acts on $V$ with $\geq 2$ Jordan blocks of size $p$.
\item $V$ is irreducible.
\item $V$ is isomorphic to a Frobenius twist of $V_G(c)$ or $V_G(c)^*$, where $p \leq c \leq 2p-2$. Furthermore, $u$ acts on $V$ with Jordan form $J_p \oplus J_{c-p+1}$.
\end{enumerate} 
\end{prop}

\begin{proof}
Let $V$ be a counterexample of minimal dimension to the claim. Then $V$ is not irreducible and $u$ acts on $V$ with $\leq 1$ Jordan block of size $p$. Note also that by Lemma \ref{jordanblockfiltration}, the element $u$ acts on any subquotient of $V$ with $\leq 1$ Jordan block of size $p$. Therefore any proper subquotient of $V$ must be as in (ii) or (iii) of the claim.

Since $V$ is not irreducible, there exists a subquotient $Q$ of $V$ which is a nonsplit extension of two irreducible $G$-modules. By Proposition \ref{nonsplitaction}, the subquotient $Q$ is isomorphic to $V_G(c)^{[l]}$ or $(V_G(c)^*)^{[l]}$ for some $l \geq 0$ and $p \leq c \leq 2p-2$. 

We are assuming that $V$ is a counterexample, so there must be a subquotient $Q'$ of $V$ which is a nonsplit extension of $Q$ and some irreducible $Z$. It is straightforward to see that such a subquotient $Q'$ must be indecomposable, and so by the minimality of $V$ we have $Q' = V$. By replacing $V$ with $V^*$ if necessary, this reduces us to the situation where $V$ is a nonsplit extension of $V_G(c)^{[l]}$ and some irreducible $L_G(d)$.

There must be a subquotient of $V$ which is a nonsplit extension of $L_G(d)$ with $L_G(c)^{[l]}$ or $L_G(2p-2-c)^{[l]}$. Therefore by Proposition \ref{nonsplitaction}, it follows that $L_G(d) = L_G(c)^{[l]}$ or $L_G(d) = L_G(2p-2-c)^{[l]}$, respectively. Thus $V \cong T_G(c)^{[l]}$ by Lemma \ref{someextlemma}. This gives us a contradiction, because $u$ acts on $T_G(c)$ with two Jordan blocks of size $p$ by Lemma \ref{lemma:tiltingSL2_jordanblocks}.\end{proof}

\begin{seur}\label{lessthanp}
Let $V$ be any representation of $G$. Suppose that $u$ acts on $V$ with all Jordan blocks of size $< p$. Then $V$ is semisimple.
\end{seur}

\begin{proof}
By Proposition \ref{fewblocksofsizep}, the only indecomposable $G$-modules on which $u$ acts with all Jordan blocks of size $< p$ are the irreducible ones.
\end{proof}

\begin{prop}\label{SL2_semisimplicity}
Let $V$ be a self-dual representation of $G$. Suppose that $u$ acts on $V$ with at most $1$ Jordan block of size $p$. Then $V$ is semisimple.
\end{prop}

\begin{proof}
Write $V = W_1 \oplus \cdots \oplus W_t$, where $W_i$ are indecomposable $G$-modules. Suppose that $W_i \not\cong W_i^*$ for some $i$. Now irreducible $G$-modules are self-dual, so $W_i$ is not irreducible and thus $u$ acts on $W_i$ with exactly $1$ Jordan block of size $p$ (propositions \ref{nonsplitaction} and \ref{fewblocksofsizep}). On the other hand, $V \cong V^* \cong W_1^* \oplus \cdots \oplus W_t^*$, so we have $W_j \cong W_i^*$ for some $i \neq j$ since the indecomposable summands are unique by the Krull-Schmidt theorem. But then $W_i \oplus W_i^*$ is a summand of $V$, which is a contradiction since $u$ acts on $W_i \oplus W_i^*$ with two Jordan blocks of size $p$. 

Thus $W_i \cong W_i^*$ for all $i$, and so by Proposition \ref{fewblocksofsizep} each $W_i$ must be irreducible, since $V_G(c) \not\cong V_G(c)^*$ for $p \leq c \leq 2p-2$.\end{proof}

We end this section with two specific results for $p = 2$, which will be used in Section \ref{section:nonGcr_orderp} to study the subgroups of $\Sp_4(K)$.

\begin{lemma}\label{lemma:sl2_T2}
Assume $p = 2$. For $V = T_G(2)$, the following properties hold:
\begin{enumerate}[\normalfont (i)]
\item $V \cong E \otimes E$.
\item $V$ has a non-degenerate $G$-invariant alternating bilinear form $\beta$, unique up to a scalar multiple.
\item The is an orthogonal decomposition $V \downarrow K[u] = J_2 \perp J_2$ with respect to any non-degenerate $G$-invariant alternating bilinear form on $V$.
\end{enumerate}
\end{lemma}

\begin{proof}
Claim (i) follows from \cite[Lemma 4]{ErdmannHenke_RINGEL}. For (ii), the claims follow from Lemma \ref{lemma:forms_quasiminuscule} (iii).

For (iii), one can proceed by explicit calculations. By (i), we can assume that $V = E \otimes E$, and by (ii) it will be enough to prove the claim for a bilinear form $\beta$ on $V$ given by the tensor product of two non-degenerate alternating bilinear forms on $E$. Then for $v = y \otimes y$, the subspace $W = \langle v, u \cdot v \rangle$ is a non-degenerate $u$-invariant subspace with $W \downarrow K[u] = J_2$, and the orthogonal decomposition $V \downarrow K[u] = J_2 \perp J_2$ is given by $V = W \perp W^\perp$.\end{proof}

\begin{lemma}\label{lemma:char2_sl2_dim4}
Assume $p = 2$. Let $V$ be a $G$-module such that $V \downarrow \fieldsymbol [u] = J_2 \oplus J_2$. Then one of the following holds:

\begin{enumerate}[\normalfont (i)]
\item $V$ is irreducible and isomorphic to $L_G(1)^{[n]} \otimes L_G(1)^{[m]}$, where $0 \leq n < m$.
\item $V \cong L_G(1)^{[n]} \oplus L_G(1)^{[m]}$, where $0 \leq n \leq m$.
\item $V \cong T_G(2)^{[n]} \cong L_G(1)^{[n]} \otimes L_G(1)^{[n]}$, where $n \geq 0$.
\end{enumerate}
\end{lemma}

\begin{proof}
If $V$ is irreducible, then the fact that $V$ has dimension $4$ implies that $V \cong L_G(1)^{[n]} \otimes L_G(1)^{[m]}$ for some $0 \leq n < m$, so $V$ is as in (i). Suppose then that $V$ is not irreducible. In this case the possible composition factors of $V$ are $L_G(0)$ and $L_G(1)^{[n]}$ for some $n \geq 0$. Consequently if $V$ is semisimple, it follows that $V$ is as in (ii) since $V \downarrow \fieldsymbol [u] = J_2 \oplus J_2$. 

Consider then the case where $V$ is not semisimple. In this case, there exists a subquotient $Q$ of $V$ which is a non-split extension between two irreducible $G$-modules. Now $Q$ has to be a proper subquotient of $V$, since there are no nonsplit extensions between $L_G(1)^{[n]}$ and $L_G(1)^{[m]}$ (Theorem \ref{sl2split}). Thus $u$ acts on $Q$ with exactly one Jordan block of size $2$, and so $Q$ must be isomorphic to $V_G(2)^{[n]}$ or $(V_G(2)^*)^{[n]}$ by Proposition \ref{fewblocksofsizep}. By replacing $V$ with $V^*$ if necessary, we may assume that $Q$ is isomorphic to $V_G(2)^{[n]}$. Since $u$ acts on $Q$ with a Jordan block of size $1$, it follows that $V$ is a nonsplit extension of $V_G(2)^{[n]}$ and $L_G(0)$. Thus $V \cong T_G(2)^{[n]}$ by Lemma \ref{someextlemma}. Finally, by Lemma \ref{lemma:sl2_T2} (i) we have $L_G(1) \otimes L_G(1) = T_G(2)$. Therefore $L_G(1)^{[n]} \otimes L_G(1)^{[n]} \cong T_G(2)^{[n]}$, which completes the proof of the lemma.\end{proof}

%% file: PROJECT_nonGcr.tex
%


In this section, we will prove our main result (Theorem \ref{theorem:main_result_of_paper}). We begin with the following basic useful observation.

\begin{lemma}\label{lemma:basicgcrdist}
Let $X < G$ be a reductive subgroup of $G$. Suppose that $X$ contains a distinguished unipotent element of $G$. Then $X$ is $G$-ir or $X$ is non-$G$-cr.
\end{lemma} 

\begin{proof}
Suppose that $X$ is not $G$-ir, i.e., that $X$ is contained in some proper parabolic subgroup $P$ of $G$. Since every Levi factor of $P$ is a centralizer of some non-trivial torus, and since $X$ contains a distinguished unipotent element, it follows that $X$ cannot be contained in any Levi factor of $P$. Hence $X$ is non-$G$-cr.\end{proof} 

As seen from the next result, for classical groups the concept of $G$-cr subgroups can be seen as a generalization of semisimplicity in representation theory. See also \cite[pg. 32-33]{LiebeckSeitzReductive}. For exceptional groups, reductive non-$G$-cr subgroups only occur in small characteristic \cite[Theorem 1]{LiebeckSeitzReductive}.

\begin{lause}[{\cite[3.2.2]{Serre}}]\label{thm:gcr_classical}
Let $G = \SL(V)$, $G = \Sp(V)$, or $G = \SO(V)$. Assume that $p > 2$ if $G = \Sp(V)$ or $G = \SO(V)$. Then a closed subgroup $H < G$ is $G$-cr if and only if $V \downarrow H$ is semisimple.
\end{lause}

In the next lemma, we will consider non-$G$-cr subgroups of type $A_1$ for simple $G$ of exceptional type in good characteristic. For the unipotent classes in $G$, we will use the labeling given by the Bala-Carter classification of unipotent conjugacy classes \cite{BalaCarter1, BalaCarter2}, which is valid in good characteristic by Pommerening's theorem \cite{Pommerening1, Pommerening2}.

\begin{lemma}\label{lemma:good_nonGcr_orderp_exceptional}
Let $G$ be a simple algebraic group of exceptional type and assume that $p$ is good for $G$. Let $u \in G$ be a unipotent element of order $p$. Then $u$ is contained in a non-$G$-cr subgroup $X < G$ of type $A_1$ precisely in the following cases:

\begin{enumerate}[\normalfont (i)]
\item $G = E_6$, $p = 5$, and $u$ is in the unipotent class $A_4$ or $A_4A_1$.
\item $G = E_7$, $p = 5$, and $u$ is in the unipotent class $A_4$, $A_4A_1$, or $A_4A_2$.
\item $G = E_7$, $p = 7$, and $u$ is in the unipotent class $A_6$.
\item $G = E_8$, $p = 7$, and $u$ is in the unipotent class $A_6$ or $A_6A_1$.
\end{enumerate}

\end{lemma}

\begin{proof}The main result of \cite{LitterickThomasARXIV} gives a complete list of non-$G$-cr subgroups $X < G$ of type $A_1$, up to $G$-conjugacy. Thus to prove our claim, it will be enough to check for each $X$ which conjugacy class of unipotent elements of order $p$ it intersects.

For each non-$G$-cr subgroup $X$, Litterick and Thomas give the $X$-module structure of the restriction of the adjoint representation of $G$, see \cite[Table 11 - Table 16]{LitterickThomasARXIV}. From their tables, we see that the adjoint representation of $G$ decomposes as an $X$-module into a direct sum of modules involving Frobenius twists, duals, and tensor products of irreducibles, Weyl modules, and tilting modules of $X$. Hence by using the decompositions in \cite[Table 11 - Table 16]{LitterickThomasARXIV}, along with Lemma \ref{symmetricpowerlemma}, Lemma \ref{lemma:tiltingSL2_jordanblocks}, and Lemma \ref{lemmatensordecomp_orderp}, we can compute the Jordan block sizes of a non-identity unipotent element $u \in X$ on the adjoint representation of $G$. Then by \cite[Theorem 2]{Lawther}, we can use the tables in \cite{Lawther} to identify the precise conjugacy class of $u$ in $G$. Doing this straightforward (but perhaps tedious) computation for each non-$G$-cr subgroup of type $A_1$ given in \cite{LitterickThomasARXIV}, one finds that they can only contain unipotent elements listed in (i) - (iv), and that all of the unipotent elements in (i) - (iv) are contained in some non-$G$-cr subgroup of type $A_1$.

We give one example of how the computation is done, all the other computations use similar methods. Let $p = 5$ and $G = E_6$. We consider an $A_1$ subgroup $X$ of $G$, which is embedded into a Levi factor of type $D_5$ via the indecomposable tilting module $T_X(8)$. According to \cite{LitterickThomasARXIV}, the subgroup $X$ is non-$G$-cr, and by \cite[Table 11]{LitterickThomasARXIV} restriction of the adjoint representation of $G$ to $X$ decomposes into a direct sum $$L_X(14) + T_X(10) + V_X(10) + V_X(10)^* + T_X(6) + L_X(4) + L_X(4) + L_X(0).$$ Let $u \in X$ be a non-identity unipotent element of $X$. We proceed to find the $\fieldsymbol[u]$-module decomposition for each of the summands. 


\begin{itemize}
\item $L_X(14)$: By Steinberg's tensor product theorem, we have $L_X(14) \cong L_X(2)^{[1]} \otimes L_X(4)$. Thus $L_X(14) \downarrow K[u] = J_3 \otimes J_5 = 3 \cdot J_5$ by Lemma \ref{lemmatensordecomp_orderp}.
\item $T_X(10)$: With Theorem \ref{thm:tiltingSL2}, one computes that $T_X(10)$ has dimension $20$, so by Lemma \ref{lemma:tiltingSL2_jordanblocks} we have $T_X(10) \downarrow \fieldsymbol [u] = 4 \cdot J_5$.
\item $V_X(10)$ and $V_X(10)^*$: For these summands, the action of $u$ has Jordan form $2 \cdot J_5 \oplus J_1$ by Lemma \ref{symmetricpowerlemma}.
\item $T_X(6)$: Here $T_X(6) \downarrow \fieldsymbol [u] = 2 \cdot J_5$ by Lemma \ref{lemma:tiltingSL2_jordanblocks}.
\item $L_X(4)$ and $L_X(0)$: Here the action of $u$ has Jordan form $J_5$ and $J_1$, respectively.
\end{itemize}

Hence $u$ acts on the adjoint representation of $G$ with Jordan form $3 \cdot J_1 \oplus 15 \cdot J_{5}$, so by \cite[Theorem 2, Table 6]{Lawther} it lies in the conjugacy class $A_4$ of $G$.\end{proof}

\begin{remark}
Let $G$ be a simple group of exceptional type and suppose that $p$ is good for $G$. Another way to phrase Lemma \ref{lemma:good_nonGcr_orderp_exceptional} is as follows: a unipotent element $u \in G$ is contained in a non-$G$-cr subgroup of type $A_1$ if and only if $A_{p-1}$ occurs in the Bala-Carter label associated with $u$. 

We omit the proof, but one can show that this is also true in the case where $G$ is simple of classical type, with a unique exception given by the case where $u \in G$ is a regular unipotent element and $G$ is simple of type $A_{p-1}$.
\end{remark}

\begin{lause}\label{thm:orderp_good_gcr}
Let $G$ be a simple algebraic group and assume that $p$ is good for $G$. Let $u \in G$ be a distinguished unipotent element of order $p$. If $X < G$ is a connected reductive subgroup of $G$ containing $u$, then $X$ is $G$-ir.
\end{lause}

\begin{proof}
Let $X$ be a connected reductive subgroup of $G$ containing $u$. We consider first the case where $X$ is simple of type $A_1$. By Lemma \ref{lemma:basicgcrdist}, it will be enough to show that $X$ is $G$-cr. If $G$ is simple of exceptional type, it is immediate from Lemma \ref{lemma:good_nonGcr_orderp_exceptional} that $X$ is $G$-cr. If $G$ is simple of classical type, then by Lemma \ref{lemma:dist_SL}, Lemma \ref{lemma:dist_odd_SP}, and Lemma \ref{lemma:dist_odd_SO}, the element $u$ acts on the natural module $V$ of $G$ with at most one Jordan block of size $p$. Now by Proposition \ref{SL2_semisimplicity} the restriction $V \downarrow X$ is semisimple, and so by Theorem \ref{thm:gcr_classical} the subgroup $X$ is $G$-cr.

Consider then the general case where $X$ is a connected reductive subgroup of $G$ containing $u$. Since $u$ is not centralized by a nontrivial torus, the same must be true for $X$, so it follows that $X$ is semisimple. Write $X = X_1 \cdots X_t$, where the $X_i$ are simple and commute pairwise. If $p \neq 3$ or if no $X_i$ is of type $G_2$, it follows from \cite[Theorem 0.1]{TestermanA1} and \cite[Theorem 5.1]{ProudSaxlTesterman} that there exists a connected simple subgroup $X' < X$ of type $A_1$ such that $u \in X'$. It follows from the previous paragraph that $X'$ is $G$-ir, so $X$ must be $G$-ir as well. 

Suppose then that $p = 3$ and that some $X_i$ is of type $G_2$. Since $p$ is good for $G$, it follows that $G$ is simple of classical type. Let $V$ be the natural module for $G$. Now $u$ has order $3$, so the largest Jordan block of size of $u$ is at most $3$. Furthermore, since $u$ is distinguished, by Lemma \ref{lemma:dist_odd_SP} and Lemma \ref{lemma:dist_odd_SO} the Jordan block sizes of $u$ are distinct and of the same parity. Hence $\dim V \leq 4$. Since every non-trivial representation of a simple algebraic group of type $G_2$ has dimension $> 4$ (see e.g. \cite{Lubeck}), it follows that no $X_i$ can be simple of type $G_2$, contradiction.\end{proof}

What remains then is to consider distinguished unipotent elements of order $p$ in bad characteristic. There are only two cases where such elements exist (type $C_2$ for $p = 2$, and type $G_2$ for $p = 3$, see proof of Theorem \ref{theorem:main_result_of_paper} below). The next proposition, which was stated in the introduction, will deal with type $C_2$ for $p = 2$. We restate it here for convenience of the reader.

\propSP*
																																											
\begin{proof}We note first that the fact that $G$ has a unique conjugacy class of distinguished unipotent elements of order $p$ follows from \cite[Proposition 6.1]{LiebeckSeitzClass}. Furthermore, by \cite[Proposition 6.1]{LiebeckSeitzClass} this conjugacy class consists precisely of those unipotent elements $u \in G$ which admit an orthogonal decomposition $V \downarrow K[u] = J_2 \perp J_2$.

Suppose that $X < G$ is connected reductive and $u \in X$. We show that either $X$ is $G$-ir, or one of (i) or (ii) holds. 

If $X$ is normalized by a maximal torus of $G$, then $X$ is $G$-cr by \cite[Proposition 3.20]{BMSGeometric}. Since $u$ is distinguished, it follows from Lemma \ref{lemma:basicgcrdist} that $X$ is $G$-ir.

Suppose then that $X$ is not normalized by any maximal torus of $G$. Since $G$ has rank $2$, it follows that $X$ is simple of type $A_1$. Now there exists a rational representation $\rho: \SL_2(K) \rightarrow \SL(V)$ such that $\rho(\SL_2(K)) = X$. If $V$ is an irreducible $X$-module, then by \cite[pg. 32-33]{LiebeckSeitzReductive} the subgroup $X$ is $G$-ir. Thus we may assume that $V$ is non-irreducible. Then by Lemma \ref{lemma:char2_sl2_dim4}, as an $\SL_2(K)$-module $V$ must be isomorphic to either

\begin{enumerate}[(1)]
\item $L_{\SL_2(K)}(1)^{[n]} \oplus L_{\SL_2(K)}(1)^{[m]}$, where $0 \leq n \leq m$; or
\item $T_{\SL_2(K)}(2)^{[n]} \cong L_{\SL_2(K)}(1)^{[n]} \otimes L_{\SL_2(K)}(1)^{[n]}$, where $n \geq 0$.
\end{enumerate}


In both cases we have $\rho = \rho' \circ F^{n}$, where $F$ is the usual Frobenius endomorphism and $\rho': \SL_2(K) \rightarrow \SL(V)$ is a rational representation. Since applying a Frobenius twist does not change the image of the representation $\rho'$, it follows that we may assume that as an $\SL_2(K)$-module, $V$ is isomorphic to either

\begin{enumerate}[(1)$'$]
\item $L_{\SL_2(K)}(1) \oplus L_{\SL_2(K)}(1)^{[n]}$, where $n \geq 0$; or
\item $T_{\SL_2(K)}(2)$.
\end{enumerate}


In case (2)$'$ we have $X$ as in case (ii) of the claim. Consider then case (1)$'$. Here if $n > 0$, it follows from \cite[pg. 32-33]{LiebeckSeitzReductive} that $X$ is $G$-ir. Suppose then that $n = 0$, so $V \cong L_{\SL_2(K)}(1) \oplus L_{\SL_2(K)}(1)$. If $V$ has an $X$-submodule $W \cong L_{\SL_2(K)}(1)$ such that $W$ is non-degenerate, then it is clear that $X$ is as in case (i) of the claim. The other possibility is that every $X$-submodule $W \cong L_{\SL_2(K)}(1)$ is totally isotropic. In this case we can find a decomposition $V = W \oplus W'$ of $X$-modules where $W \cong L_{\SL_2(K)}(1) \cong W'$ and $W$, $W'$ are totally isotropic. But then $X$ is a Levi factor of type $A_1$ (short root), which contradicts the assumption that $u$ is distinguished.


We consider the existence and uniqueness claims for the subgroups $X$ in (i) and (ii). We begin by considering subgroups $X$ in (i). First, note that the $\SL_2(K)$-module $L_{\SL_2(K)}(1)$ has a non-degenerate $\SL_2(K)$-invariant alternating bilinear form. Therefore it is clear that we can find a representation $\rho: \SL_2(K) \rightarrow G$ with $V = V_1 \perp V_2$ such that $V_1 \cong L_{\SL_2(K)}(1) \cong V_2$, so we can choose $X = \rho(\SL_2(K))$. The uniqueness of such an $X$ up to $G$-conjugacy follows easily from the fact that any two orthogonal direct sum decompositions $V_1 \perp V_2$ and $W_1 \perp W_2$ are conjugate under the action of $G$.


Next, note that since a non-identity unipotent element $u \in X$ acts on the module $L_{\SL_2(K)}(1)$ with a single Jordan block of size $2$, it is clear that $V \downarrow \fieldsymbol [u] = J_2 \perp J_2$. For subgroups $X$ in (i), the fact that $X$ is contained in a proper parabolic subgroup of $G$ follows from \cite[pg. 32-33]{LiebeckSeitzReductive}. This completes the proof of the claims for (i).

For the subgroup $X$ in (ii), note that by Lemma \ref{lemma:sl2_T2} (ii) the tilting $\SL_2(K)$-module $T_{\SL_2(K)}(2)$ has a non-degenerate $\SL_2(K)$-invariant alternating bilinear form. Therefore there exists a representation $\rho: \SL_2(K) \rightarrow G$, with $V \downarrow X \cong T_{\SL_2(K)}(2)$ for $X = \rho(\SL_2(K))$. Also by Lemma \ref{lemma:sl2_T2} (ii), such an $X$ is unique up to $G$-conjugacy. Next note that for a non-identity unipotent element $u \in X$, we have $V \downarrow \fieldsymbol[u] = J_2 \perp J_2$ by Lemma \ref{lemma:sl2_T2} (iii). Finally, since $T_{\SL_2(K)}(2)$ is indecomposable and non-irreducible, it follows that $X$ is contained in a proper parabolic subgroup of $G$ \cite[pg. 32-33]{LiebeckSeitzReductive}. This completes the proof of the claims for (ii).\end{proof}


We are now ready to prove our main result, which we also restate for convenience.



\mainresult*

\begin{proof}
Suppose that $X$ is contained in a proper parabolic subgroup of $G$. Since $u$ is a distinguished unipotent element, $X$ must be non-$G$-cr (Lemma \ref{lemma:basicgcrdist}). Then by Theorem \ref{thm:orderp_good_gcr} we have that $p$ is bad for $G$. 

Consider first the case where $G$ is simple of exceptional type. Looking at the tables in \cite{Lawther}, the fact that $p$ is bad for $G$ and the fact that $u$ is a distinguished unipotent element of order $p$ implies that $G = G_2$ and $p = 3$. However, in this case it follows from \cite[Corollary 2]{StewartG2} that every connected reductive subgroup of $G$ is $G$-cr, contradicting the fact that $X$ is non-$G$-cr.

Suppose then that $G$ is simple of classical type. Since $p$ is bad for $G$, we have $p = 2$ and $G$ is simple of type $B_l$ ($l \geq 3$), $C_l$ ($l \geq 2$), or $D_l$ ($l \geq 4$). Let $V = L_G(\omega_1)$ be the natural irreducible representation of $G$. Now $u$ is a unipotent element of order $2$, so $u$ acts on $V$ with largest Jordan block of size at most $2$. Furthermore, since $u$ is a distinguished unipotent element, by Lemma \ref{lemma:dist_even_SP} we have $V \downarrow K[u] = J_2 \perp J_2$. Hence $\dim V = 4$ and $G$ is simple of type $C_2$, and in this case the claim follows from Proposition \ref{prop:resultSP4_intro}.\end{proof}


%% file: PROJECT_examples_v2.tex
%

In this final section, we give further examples of connected reductive subgroups which contain a distinguished unipotent element and which are contained in a proper parabolic subgroup. All of the examples that we give here arise from indecomposable tilting modules with quasi-minuscule highest weight, which were studied in Section \ref{section:quasiminuscule}. Other examples also exist, see for example Proposition \ref{prop:resultSP4_intro} (i). However, we believe that there are not too many examples, and aim to classify all of them in future work.

We begin with the following example, which is the only one of which we are aware in characteristic $p \neq 2$.

\begin{prop}\label{prop:typeF4_nonGcr}
Assume that $p = 3$ and let $G = \SO(V)$ with $\dim V = 27$. Then:

\begin{enumerate}[\normalfont (i)]
\item Up to $G$-conjugacy there exists a unique $X < G$ simple of type $F_4$ such that $V \downarrow X \cong T_X(\omega_4)$;
\item Such an $X$ is contained in a proper parabolic subgroup of $G$;
\item A regular unipotent element $u \in X$ satisfies $V \downarrow K[u] = J_3 \oplus J_9 \oplus J_{15}$, and thus is a distinguished unipotent element of $G$ (Lemma \ref{lemma:dist_odd_SO}).
\end{enumerate}
\end{prop}

\begin{proof}
Let $X$ be a simple algebraic group of type $F_4$. First note that the weight $\omega_4$ is quasi-minuscule and $V_X(\omega_4) = L_X(\omega_4) | L_X(0)$ (see Table \ref{table:quasiminuscule}) has dimension $26$. Therefore by Lemma \ref{lemma:forms_quasiminuscule} (i), the indecomposable tilting module $T_X(\omega_4)$ is $27$-dimensional, uniserial, and $T_X(\omega_4) = L_X(0) | L_X(\omega_4) | L_X(0)$. Furthermore, we have a non-degenerate $G$-invariant symmetric bilinear form on $T_X(\omega_4)$ by Lemma \ref{lemma:forms_quasiminuscule} (ii). This shows the existence part in claim (i). Since $p \neq 2$, uniqueness up to conjugacy is a consequence of a general result in representation theory, see \cite[Korollar 3.5]{Quebbemann} or alternatively \cite[Satz 3.12]{WillemsThesis}. This establishes (i). Claim (ii) follows from the fact that $T_X(\omega_4)$ is an indecomposable and non-irreducible $X$-module \cite[pg. 32-33]{LiebeckSeitzReductive}.

To establish (iii), we use the usual embedding of $X$ into a simply connected simple algebraic group $Y$ of type $E_6$. That is, we consider $X$ as the centralizer of the involutory graph automorphism of $Y$ induced by the nontrivial automorphism of the Dynkin diagram of $Y$. 

Since in type $E_6$ the fundamental highest weight $\omega_1$ is minuscule, the Weyl module $V_Y(\omega_1)$ is irreducible and thus we have $V_Y(\omega_1) = T_Y(\omega_1)$. According to \cite[Theorem 20]{vanderKallen}, the restriction of every tilting module of $Y$ to $X$ is a tilting module for $X$. In particular, the restriction $V_Y(\omega_1) \downarrow X$ is tilting for $X$. For the character of this restriction, we have $$\ch V_Y(\omega_1) \downarrow X = \ch V_X(\omega_4) + \ch V_X(0)$$ by \cite[Table 8.7]{LiebeckSeitzReductive}. Since $V_Y(\omega_1) \downarrow X$ is tilting, we conclude that $V_Y(\omega_1) \downarrow X \cong T_X(\omega_4)$.

Let $u \in X$ be a regular unipotent element. It follows from \cite[Theorem A]{SaxlSeitz} that $u$ is also a regular unipotent element of $Y$. Then according to \cite[Table 5]{Lawther}, we have $V_Y(\omega_1) \downarrow K[u] = J_3 \oplus J_9 \oplus J_{15}$, hence $T_X(\omega_4) \downarrow K[u] = J_3 \oplus J_9 \oplus J_{15}$.\end{proof}

We now give more examples of similar nature in characteristic two.

%
%

\begin{prop}\label{prop:nonGcr_G2_C4}
Assume that $p = 2$ and let $G = \Sp(V)$ with $\dim V = 8$. Then:

\begin{enumerate}[\normalfont (i)]
\item Up to $G$-conjugacy there exists a unique $X < G$ simple of type $G_2$ such that $V \downarrow X \cong T_X(\omega_1)$;
\item Such an $X$ is contained in a proper parabolic subgroup of $G$;
\item A regular unipotent element $u \in X$ satisfies $V \downarrow K[u] = J_2 \perp J_6$, and thus is a distinguished unipotent element of $G$ (Lemma \ref{lemma:dist_even_SP}).
\end{enumerate}
\end{prop}

\begin{prop}\label{prop:nonGcr_E7}
Assume that $p = 2$ and let $G = \Sp(V)$ with $\dim V = 134$. Then:

\begin{enumerate}[\normalfont (i)]
\item Up to $G$-conjugacy there exists a unique $X < G$ simple of type $E_7$ such that $V \downarrow X \cong T_X(\omega_1)$;
\item Such an $X$ is contained in a proper parabolic subgroup of $G$;
\item A regular unipotent element $u \in X$ satisfies $V \downarrow K[u] = J_2 \perp J_8 \perp J_{10} \perp J_{16} \perp J_{18} \perp J_{22} \perp J_{26} \perp J_{32}$, and thus is a distinguished unipotent element of $G$ (Lemma \ref{lemma:dist_even_SP}).
\end{enumerate}
\end{prop}

\begin{proof}[{Proof of Propositions \ref{prop:nonGcr_G2_C4} and \ref{prop:nonGcr_E7}.}]Let $X$ be simple of type $G_2$ or $E_7$. The fundamental dominant weight $\omega_1$ is quasi-minuscule, and the corresponding Weyl module has structure $V_X(\omega_1) = L_X(\omega_1) | L_X(0)$ (Table \ref{table:quasiminuscule}). Thus (i) follows from Lemma \ref{lemma:forms_quasiminuscule} (i) and (iii), along with the fact that $\dim V_X(\omega_1) = 7$ if $X$ has type $G_2$ and $\dim V_X(\omega_1) = 133$ if $X$ has type $E_7$. Claim (ii) follows from the fact that $T_X(\omega_1)$ is indecomposable and non-irreducible \cite[pg. 32-33]{LiebeckSeitzReductive}.

For the proof of (iii), write $V = T_X(\omega_1)$ and let $\beta$ be a non-degenerate $X$-invariant bilinear form on $V$, given by Lemma \ref{lemma:forms_quasiminuscule} (iii). Since $V$ is an indecomposable tilting $X$-module, it follows from \cite[Corollary 7.2]{GaribaldiNakano} that there exists an $X$-invariant quadratic form $Q$ on $V$ such that $\beta$ is the polarization of $Q$. Let $v \in V$ be a vector generating the unique $1$-dimensional $X$-submodule of $V$. Then $\langle v \rangle^\perp / \langle v \rangle \cong L_X(\omega_1)$ by Lemma \ref{lemma:forms_quasiminuscule} (i). Note that $L_X(\omega_1)$ does not admit a non-degenerate $X$-invariant quadratic form (see e.g. \cite[Proposition 6.1]{Korhonen} and \cite[Theorem 3.4 (a)]{GowWillems1}), so $Q(v) \neq 0$.


Let $u \in X$ be a regular unipotent element. We describe the action of $u$ on $\langle v \rangle^\perp / \langle v \rangle \cong L_X(\omega_1)$, which together with Lemma \ref{lemma:typeBtoDstabilizer} will complete the proof of (iii). If $X$ has type $G_2$, it follows from \cite[Theorem 1.9]{Suprunenko} that $L_X(\omega_1) \downarrow K[u] = J_6$. For $X$ of type $E_7$, note that $\langle v \rangle^\perp = V_X(\omega_1)$. We have $$\langle v \rangle^\perp \downarrow K[u] = J_1 \oplus J_8 \oplus J_{10} \oplus J_{16} \oplus J_{18} \oplus J_{22} \oplus J_{26} \oplus J_{32}$$ by \cite[Table 8]{Lawther}, because $V_X(\omega_1)$ is the adjoint representation or its dual. As noted in the proof of Lemma \ref{lemma:typeBtoDstabilizer}, we can find an orthogonal decomposition $\langle v \rangle^\perp \downarrow K[u] = W \perp \langle v \rangle$ with $W \downarrow K[u] = J_8 \oplus J_{10} \oplus J_{16} \oplus J_{18} \oplus J_{22} \oplus J_{26}$, which combined with Lemma \ref{lemma:dist_even_SP_noreps} gives $$\langle v \rangle^\perp / \langle v \rangle \downarrow K[u] = J_8 \perp J_{10} \perp J_{16} \perp J_{18} \perp J_{22} \perp J_{26}$$ and completes the proof.\end{proof}


%

\begin{remark}\label{remark:infinite}An infinite family of examples, similar to Proposition \ref{prop:nonGcr_G2_C4} and Proposition \ref{prop:nonGcr_E7}, is given by the following. Assume that $p = 2$ and consider the adjoint simple algebraic group $X$ of type $B_{l-1}$, where $l \geq 2$. Let $G = \Sp(V, \beta)$, where $\dim V = 2l$ and $\beta$ is a non-degenerate alternating bilinear form on $V$. Consider $Y = \SO(V, Q) < G$ where $Q$ is a non-degenerate quadratic form on $V$ that polarizes to $\beta$. Then we can view $X = \operatorname{Stab}_Y(v)$, where $v \in V$ is a non-zero vector such that $Q(v) \neq 0$.

In this setting, we have $V \downarrow X \cong T_X(\omega_1)$, and $X$ is contained in a proper parabolic subgroup since it stabilizes the totally isotropic subspace $\langle v \rangle$. Furthermore, we can see that $X$ contains distinguished unipotent elements of $G$. Let $u \in X$ be a unipotent element such that with respect to the bilinear form induced on $\langle v \rangle^\perp / \langle v \rangle$, we have $$\langle v \rangle^\perp / \langle v \rangle \downarrow K[u] = J_{d_1} \perp \cdots \perp J_{d_t},$$ where $t$ is odd, each $d_i$ is even, each $J_{d_i}$ occurs with multiplicity at most two, and $J_2$ occurs with multiplicity at most one. Then $$V \downarrow K[u] = J_2 \perp J_{d_1} \perp \cdots \perp J_{d_t}$$ by Lemma \ref{lemma:typeBtoDstabilizer}, so $u$ is a distinguished unipotent element of $G$ by Lemma \ref{lemma:dist_even_SP}. For example, for a regular unipotent element $u$ of $X$, we have $V \downarrow K[u] = J_2 \perp J_{2l-2}$.\end{remark}

\begin{remark}\label{remark:furthertodo}
We describe one more example arising from an indecomposable tilting module with quasi-minuscule highest weight. Assume that $p = 2$. One can show (see Remark \ref{remark:quasi_typeD_omitted}) that in $G = \Sp(V)$ for $\dim V = 68$, up to $G$-conjugacy there exists a unique $X < G$ simple of type $D_6$ such that $V \downarrow X \cong T_X(\omega_2) = L_X(0)^2 | L_X(\omega_2) | L_X(0)^2$. We omit the proof from this paper, but one can show that $X$ is contained in a proper parabolic subgroup of $G$ and that a regular unipotent element $u \in X$ satisfies $$V \downarrow K[u] = J_2 \perp J_2 \perp J_6 \perp J_8 \perp J_{10} \perp J_{10} \perp J_{14} \perp J_{16}$$ with respect to the alternating bilinear form of $G$. Note that then $u$ is a distinguished unipotent element of $G$ by Lemma \ref{lemma:dist_even_SP}.
\end{remark}


%% file: PROJECT_acknowledgements.tex
%

The results of this paper were obtained during my doctoral studies, which was supported by a grant from the Swiss National Science Foundation (grant number $200021 \_ 146223$). I am very thankful to my advisor, Prof. Donna Testerman, for her guidance and for her many useful comments on the content this paper.
